\documentclass[12pt]{amsart}

\usepackage{amssymb,amsmath,amsthm,enumerate,graphicx,stmaryrd,tikz}
\usetikzlibrary{arrows}

\newcommand{\ox}{\mathcal{O}_X}

\newcommand{\st}{\mathrm{st}}
\newcommand{\semi}{\mathrm{ss}}
\newcommand{\un}{\mathrm{un}}
\newcommand{\NR}[1]{N^1(#1)_{\mathbb{R}}}
\newcommand{\NQ}[1]{N^1(#1)_{\mathbb{Q}}}
\newcommand{\PicQ}[1]{\Pic(#1)_{\mathbb{Q}}}

\def\<{\langle}
\def\>{\rangle}
\def\AA{\mathcal{A}}

\def\CC{\mathbb{C}}

\newcommand{\LL}{\mathbf{L}}
\newcommand{\mm}{\mathbf{m}}
\def\NN{\mathbb{N}}
\def\OO{\mathcal{O}}
\def\PP{\mathbb{P}}
\def\QQ{\mathbb{Q}}

\def\ZZ{\mathbb{Z}}
\newcommand{\xx}{\mathbf{x}}

\DeclareMathOperator{\Cox}{Cox}

\DeclareMathOperator{\Hom}{Hom}

\DeclareMathOperator{\Nef}{Nef}
\DeclareMathOperator{\NS}{NS}
\DeclareMathOperator{\Pic}{Pic}
\DeclareMathOperator{\Proj}{Proj}
\DeclareMathOperator{\Spec}{Spec}

\theoremstyle{plain}
\newtheorem{thm}{Theorem}[section]
\newtheorem{prop}[thm]{Proposition}
\newtheorem{lemma}[thm]{Lemma}

\newtheorem{theorem}{Theorem}
\newtheorem{corollary}[theorem]{Corollary}
\newtheorem{proposition}[theorem]{Proposition}

\theoremstyle{definition}
\newtheorem{defn}[thm]{Definition}
\newtheorem{rmk}[thm]{Remark}
\newtheorem{notation}[thm]{Notation}

\newtheorem{definition}[theorem]{Definition}
\newtheorem{remark}[theorem]{Remark}
\newtheorem{example}[theorem]{Example}

\begin{document}

\mbox{}
\vspace{-1.1ex}
\title{A Lefschetz hyperplane theorem for Mori dream spaces}
\author{Shin-Yao Jow}
\date{27 November 2009}

\begin{abstract}
 Let $X$ be a smooth Mori dream space of dimension $\ge 4$. We show that, if $X$
 satisfies a suitable GIT condition which we call \emph{small unstable locus}, then
 every smooth ample divisor $Y$ of $X$ is also a Mori dream space.
 Moreover, the restriction map identifies the N\'{e}ron-Severi spaces of $X$ and $Y$,
 and under this identification every Mori chamber of $Y$ is a union of some Mori
 chambers of $X$, and the nef cone of $Y$ is the same as the nef cone of $X$. This
 Lefschetz-type theorem enables one to construct many examples of Mori dream spaces by
 taking ``Mori dream hypersurfaces'' of an ambient Mori dream space, provided that it
 satisfies the GIT condition. To facilitate this, we then show that the GIT condition is
 stable under taking products and taking the projective bundle of the direct sum of at
 least three line bundles, and in the case when $X$ is toric, we show that the condition
 is equivalent to the fan of $X$ being $2$-neighborly.
\end{abstract}


\maketitle

\section*{Introduction}

The main purpose of this paper is to prove an analogue of the Lefschetz hyperplane theorem for
Mori dream spaces.

Let $X$ be a smooth complex projective variety, and let $N^1(X)$ be the group of
numerical equivalence classes of line bundles on $X$. Recall from \cite{HK} that $X$
is called a Mori dream space if $\PicQ{X}=\NQ{X}$ (equivalently $H^1(X,\ox)=0$),
and $X$ has a finitely generated Cox ring (Definition~\ref{d:Cox}).
As the name might suggest, Mori dream spaces are very special
varieties on which Mori theory works extremely well (see the nice survey article of
Hu \cite{Hu}). On the other hand,
not many classes of examples of them are known. It has been understood for a while that toric varieties
are Mori dream spaces; indeed their Cox rings are
polynomial rings, Cox's homogeneous coordinate rings \cite{Cox}. Besides that, it
was only proved very recently, in the spectacular paper of \cite{BCHM}, that
(log) Fano varieties are also Mori dream spaces. The Cox rings of certain
Mori dream spaces have been the focus of much study: see, for example,
\cite{BP}, \cite{STV}, \cite{SS}, \cite{CT}, \cite{Cas}.

The most prominent feature of a Mori dream space discovered in \cite{HK} is the existence of
a polyhedral chamber decomposition of its pseudo-effective cone; these chambers are known
as the \emph{Mori chambers}. Specifically if $L$ is a line
bundle on a Mori dream space $X$, then its section ring \[
 R(X,L):=\bigoplus_{n\in \NN} H^0(X,L^{\otimes n}) \]
is finitely generated. Thus the rational map defined by the linear series $|L^{\otimes n}|$ \[
 \phi_{|L^{\otimes n}|}\colon X \dashrightarrow \PP H^0(X,L^{\otimes n}) \]
stabilizes to some rational map \[
 \phi_{L}\colon X \dashrightarrow \Proj R(X,L) \]
for all large and sufficiently divisible $n$. Two line bundles $L_1$ and
$L_2$ are said to be \emph{Mori equivalent} if $\phi_{L_1}=\phi_{L_2}$. This equivalence relation naturally
extends to $\PicQ{X}$, and a Mori chamber is just the closure of an equivalence class in $\NR{X}$
which has a nonempty interior. It was shown in \cite{HK} that these Mori chambers are
polyhedral and in one-to-one correspondence with birational contractions of $X$ having
$\QQ$-factorial image.

In this paper, we first define the notion of a \emph{Mori dream hypersurface} of
a Mori dream space. Since the chamber structure plays such a key role in the geometry of
a Mori dream space, we propose that what deserved to be called a Mori dream hypersurface
should not only be a Mori dream space itself, but should also respect the chamber
structure in the following sense:

\begin{definition} \label{d:1}
 Let $X$ be a Mori dream space. A hypersurface $Y\subset X$ is called a \emph{Mori dream
 hypersurface} if it satisfies the following three requirements:
 \begin{enumerate}[(i)]
  \item $Y$ is a Mori dream space;
  \item The restriction map determines an isomorphism between $\NR{X}$ and $\NR{Y}$;
  \item After identifying $\NR{X}$ and $\NR{Y}$ via the restriction map, each Mori chamber
  of $Y$ is a union of some Mori chambers of $X$.
 \end{enumerate}
\end{definition}

Note that the second requirement in the above definition is satisfied for any smooth
projective variety $X$ of dimension $\ge 4$ and $Y\subset X$ a smooth ample divisor,
thanks to the Lefschetz hyperplane theorem. On the other hand, if $X=\PP^1\times
\PP^{n-1}$, then ample divisors $Y\subset X$ are generally not Mori dream hypersurfaces. This leads to the
question of finding suitable conditions under which a ``Lefschetz-type'' theorem would
hold \emph{in the category of Mori dream spaces}. We will give one such condition in this
paper. Before stating our condition, we give some corollaries:

\begin{corollary} \label{c:2}
 Let $X$ be a smooth projective variety of dimension $\ge 4$. Suppose $X$ is a product of some
 Mori dream spaces, each having dimension $\ge 2$ and Picard number one. Then $X$ is a
 Mori dream space, and every smooth ample divisor $Y\subset X$ is a Mori dream hypersurface;
 moreover, the restriction map identifies $\Nef(X)$, the nef cone of\/ $X$,
 with $\Nef(Y)$, the nef cone of\/ $Y$.
\end{corollary}

\begin{corollary} \label{c:3}
 Let $X$ be a smooth projective toric variety of dimension $\ge 4$ associated to a fan $\Delta$.
 Suppose that for any two rays in $\Delta$, the two-dimensional convex cone they span is also in
 $\Delta$. Then $X$ is a Mori dream space in which every smooth ample divisor $Y$ is a Mori dream
 hypersurface, and the restriction map identifies $\Nef(X)$ with $\Nef(Y)$.
\end{corollary}

In fact more examples satisfying the conclusion of the above two corollaries can be obtained
by a suitable projective bundle construction: see Proposition~\ref{p:4}.

\begin{example}
 The simplest example of a space $X$ as in Corollary~\ref{c:2} is a product of general complete
 intersections in projective spaces. The simplest example of a space $X$ in
 Corollary~\ref{c:3} other than $\PP^n$ is the blowup of $\PP^n$ along a linear subspace
 $\PP^m$ for $0<m<n-2$.
\end{example}

\begin{remark}
 In Corollary~\ref{c:2}, the part of the result about the preservation of nef cones has
 previously been obtained by Hassett-Lin-Wang \cite[Theorem~4.1]{HLW}, which they called
 ``the weak Lefschetz principle for ample cones''. See also the
 results of Koll\'{a}r \cite[Appendix]{Bor} and Wi\'sniewski \cite[Theorem~2.1]{Wis}.
 In the category of
 Mori dream spaces, however, our result applies to more spaces, such as those in
 Corollary~\ref{c:3}, which are not covered by the results in \cite{HLW}.
\end{remark}

To explain the condition lying behind the above corollaries which allows a Mori
dream space to enjoy this Lefschetz-type property for its nef cone and ample divisors, we
need to recall the GIT construction of a Mori dream space \cite[Proposition~2.9]{HK},
which says roughly that every Mori dream space $X$ is naturally a GIT quotient of an
affine variety under an algebraic torus action. More specifically, let $V=\Spec R$ where
$R$ is a Cox ring of $X$. Since $R$ is graded by a lattice $N$ in the N\'{e}ron-Severi
space of $X$, the algebraic torus $T=\Hom(N,\CC^*)$ naturally acts on the affine variety
$V$. Let $\chi\in N$ be a character of $T$ which corresponds to an ample class in the
N\'{e}ron-Severi space of $X$. Then Hu and Keel showed that $X= V\sslash_{\!\chi} T$,
the GIT quotient constructed with respect to the trivial line bundle on $V$ endowed with a
$T$-linearization by $\chi$. Moreover, this GIT quotient is a good geometric quotient, and the
unstable locus $V^\un_\chi$ always has codimension $\geq 2$ in $V$. These considerations
suggest the following theorem, which we will prove in Section~\ref{s:ThmPf}:

\begin{theorem} \label{t:main}
 Let $X$ be a smooth Mori dream space of dimension $\ge 4$, and let\/ $V$, $T$, and\/ $\chi$ be as
 above. Assume further that the following condition~\eqref{e:*} is satisfied:
 \begin{equation} \label{e:*}
  \text{The unstable locus $V^\un_\chi$ has codimension $\geq 3$ in $V$.} \tag{$\ast$}
 \end{equation}
 Then every smooth ample divisor $Y\subset X$ is a Mori dream hypersurface, and the restriction
 map identifies $\Nef(X)$ with $\Nef(Y)$.
\end{theorem}

\begin{definition}
 We will say that a Mori dream space $X$ has \emph{small unstable locus} if the
 condition~\eqref{e:*} above is satisfied.
\end{definition}

From this theorem, Corollary~\ref{c:2} and \ref{c:3} follow immediately once
the following Proposition~\ref{p:4} and \ref{p:5} are established in Section~\ref{s:PropPf}:

\begin{proposition} \label{p:4}
 Let $X$, $X_1$ and $X_2$ be Mori dream spaces.
 \begin{enumerate}[\upshape (a)]
  \item If\/ $X$ has dimension at least two and Picard number equal to one,
   then $X$ has small unstable locus.
  \item If\/ $X_1$ and $X_2$ both have small unstable locus, then
   $X_1\times X_2$ is a Mori dream space which has small unstable locus.
  \item Suppose that\/ $X$ has small unstable locus. Let $L_1,\ldots,L_k$ be line bundles on $X$,
   $k\geq 3$. Then the projective bundle $\PP(\bigoplus_{i=1}^k L_i^{\otimes m})$
   is a Mori dream space having small unstable locus for all sufficiently divisible
   integers $m$.
 \end{enumerate}
\end{proposition}

\begin{definition}
 A fan $\Delta$ is called \emph{$m$-neighborly} if for any $m$ rays in $\Delta$,
 the convex cone they span is also in $\Delta$.
\end{definition}

\begin{proposition}  \label{p:5}
 Let $X$ be a simplicial projective toric variety associated to a fan $\Delta$. Then $X$
 has small unstable locus if and only if $\Delta$ is $2$-neighborly.
\end{proposition}

\begin{remark}
 Fans which are $m$-neighborly and give rise to complete smooth toric varieties have been the subject
 of interest in a couple of papers by Kleinschmidt, Sturmfels and others (\cite{GKS},
 \cite{KSS}). Our proof of Proposition~\ref{p:5} indeed shows that the fan $\Delta$ is
 $m$-neighborly if and only if in Cox's GIT description of the corresponding toric
 variety $X$ \cite{Cox}, the unstable locus has codimension at least $m+1$. This
 reveals that the neighborliness property of the fan, which is of a
 combinatorial nature, has a nice GIT interpretation on the corresponding variety side.
\end{remark}

\begin{remark}
 We point out that using Proposition~\ref{p:4}~(c), one can construct Mori dream spaces
 which have small unstable locus and also possess a nontrivial small $\QQ$-factorial modification
 \cite[Definition~1.8]{HK}. For example let $Z$ be the blowup of $\mathbb{P}^4$ along a line, let $L_1$ be the line bundle on $Z$ corresponding to the exceptional divisor, and let $L_2$ and $L_3$ both be the pullback of $\OO_{\PP^4}(1)$ to $Z$. Then
 Proposition~\ref{p:4}~(c) says that
 $X=\mathbb{P}(\bigoplus_{i=1}^3 L_i^{\otimes m})$ is a Mori dream space which has small
 unstable locus if $m$ is sufficiently divisible. To see that $X$
 has a nontrivial small $\QQ$-factorial modification, note that the stable base locus 
 of $\mathcal{O}_X(1)$ has no divisorial component (in fact is has codimension $3$). So the moving cone of $X$ is strictly larger than $\Nef(X)$, hence $X$ must have a small $\QQ$-factorial modification other than itself \cite[Proposition~1.11]{HK}.
\end{remark}

 Finally, we remark that it would be interesting to clarify the relation between Wi\'sniewski's
 \cite[Theorem~2.1]{Wis} and our Theorem~\ref{t:main}.

\subsection*{Acknowledgements}
The author would like to thank Yi Hu, Se\'an Keel, Robert Lazarsfeld, and Mircea Musta\c{t}\v{a} for valuable
discussions and suggestions.

\section{Mori dream space as a GIT quotient} \label{s:GIT}

In this section we collect some results and set some notations which will be used later,
centering around the idea of representing a Mori dream space as a GIT quotient. We also
show in Proposition~\ref{p:normal} that every Mori dream space has a normal Cox ring (cf. \cite[Corollary~1.2]{EKW}).
We point out that the important Theorem~\ref{t:Prop 2.9} and Theorem~\ref{t:Thm 2.3} are taken
from \cite{HK}.

\begin{notation} \label{n:LL}
  For an $r$-tuple of line bundles
    $\LL=(L_1,\ldots,L_r)$ on a projective variety $X$ and an $r$-tuple of integers
    $\mathbf{m}=(m_1,\ldots,m_r)$, we let \[
    \LL^\mathbf{m}:=L_1^{\otimes m_1}\otimes L_2^{\otimes m_2}\otimes\cdots\otimes L_r^{\otimes
     m_r}. \]
  Also we let $N^1(X,\LL)\subset N^1(X)$ be the subgroup generated by $[\LL^\mm]$, the numerical
  class of $\LL^\mm$, for all $\mm\in\ZZ^r$, and we define $T_{\LL}$ to be the algebraic torus
  whose character group $\chi(T_{\LL})$ is $N^1(X,\LL)$: \[
     T_{\LL} := \Hom(N^1(X,\LL),\CC^*). \]
\end{notation}

\begin{defn} \label{d:Cox}
    Let $X$ be a projective variety such that $\PicQ{X}=\NQ{X}$.
    By a \emph{Cox ring} for $X$ we mean the ring \[
      \Cox(X,\LL):=\bigoplus_{\mm\in\ZZ^r}H^0(X,\LL^{\mm}) \]
    where $\LL=(L_1,\ldots,L_r)$ are line bundles which form a basis of $\PicQ{X}$.
    Note that the natural $N^1(X,\LL)$-grading on $\Cox(X,\LL)$ corresponds to a
    $T_{\LL}$-action on $\Spec \Cox(X,\LL)$.
\end{defn}

\begin{rmk}
  Although the definition of $\Cox(X,\LL)$ depends on a choice of basis $\LL$, whether or
  not it is \emph{finitely generated} is independent of this choice, due to the following
  well-known fact:
\end{rmk}

\begin{lemma} \label{l:R}
 Let $R$ be a $\ZZ^r$-graded commutative ring with identity. For any $\mm\in\ZZ^r$, we
 denote the subset of $R$ consisting of all degree-$\mm$ homogeneous elements and $0$ as $R_{\mm}$,
 and we define  \[
   R^{(\mm)}:=\bigoplus_{\mathbf{a}\in\ZZ^r} R_{(a_1m_1,\ldots,a_rm_r)}. \]
 If $R$ is an integral domain, and the subring $R_0:=R_{(0,\ldots,0)}$ is Noetherian,
 then the following are equivalent:
  \begin{enumerate}[\upshape (a)]
   \item $R$ is a finitely generated $R_0$-algebra;
   \item There exists an $\mm\in\ZZ_{>0}^r$ such that $R^{(\mm)}$ is a finitely generated
   $R_0$-algebra;
   \item For any $\mm\in\ZZ_{>0}^r$, $R^{(\mm)}$ is a finitely generated
   $R_0$-algebra.
  \end{enumerate}
 Moreover, when this is the case, then $R$ is a finitely generated $R^{(\mm)}$-module for
 any $\mm\in\ZZ_{>0}^r$.
\end{lemma}

\begin{proof}
 (b)$\Rightarrow$(a): As an $R^{(\mm)}$-module, $R$ is the direct sum of all modules of the form $R^{(\mm)+\mathbf{b}}$ where $\mathbf{b}\in \ZZ^{r}$ and $0\le b_i < m_i$ for all $i$. Hence it suffices to show that each of these is a finitely generated $R^{(\mm)}$-module. If $R^{(\mm)+\mathbf{b}}=0$ then this is trivial. Otherwise pick a nonzero homogeneous element $x\in R^{(\mm)+\mathbf{b}}$. Since $R$ is a domain, $R^{(\mm)+\mathbf{b}}$ is isomorphic to $x^{m_1\cdots m_r -1}R^{(\mm)+\mathbf{b}}$ as $R^{(\mm)}$-module, and since $x^{m_1\cdots m_r -1}R^{(\mm)+\mathbf{b}}\subset R^{(\mm)}$, it is finitely generated because $R^{(\mm)}$ is a Noetherian ring by Hilbert basis theorem.
 
 (a)$\Rightarrow$(c): Assuming (a), we will show that there exists $\boldsymbol{\ell}\in \ZZ_{>0}^r$ such that for all $\mm\in \ZZ_{>0}^r$, $R^{(m_1\ell_1,\ldots,m_r\ell_r)}$ is a finitely generated $R_0$-algebra. Then together with (b)$\Rightarrow$(a) we already proved, this implies (c). Suppose $R=R_0[x_1,\ldots,x_k]$, where $x_i$ is homogeneous of degree $(d_{i1},\ldots,d_{ir})\in \ZZ^r$. We will choose the $\boldsymbol{\ell}\in \ZZ_{>0}^r$ whose $j$th entry is given by 
 \begin{align*}
  \ell_j={}&\text{the positive least common multiple of the nonzero numbers in $\{d_{1j}, \ldots,d_{kj}\}$,} \\
         &\text{or $1$ if $d_{1j}=d_{2j}=\cdots =d_{kj}=0$.}
 \end{align*}
 To see that this choice works, we define, for each $i\in \{1,\ldots,k\}$, the following positive integer $h_i$: \[
  h_i=\prod_{1\le j\le r,d_{ij}\ne 0}\frac{m_j \ell_j}{|d_{ij}|}. \]
Then  \[
  R_0[x_1^{h_1},\ldots,x_k^{h_k}]\subset R^{(m_1\ell_1,\ldots,m_r\ell_r)} \subset R_0[x_1,\ldots,x_k]=R.  \]
  Since $R_0[x_1^{h_1},\ldots,x_k^{h_k}]$ is a Noetherian ring and $R_0[x_1,\ldots,x_k]$ is obviously a finitely generated module over it, the submodule $R^{(m_1\ell_1,\ldots,m_r\ell_r)}$ is thus finitely generated as well.
\end{proof}

\begin{defn}
  We will call $X$ a \emph{Mori dream space} if $X$ is a normal $\QQ$-factorial projective
    variety with $\PicQ{X}=\NQ{X}$ and a finitely generated Cox ring.
\end{defn}

\begin{rmk} \label{r:Pic=N}
 \begin{enumerate}[(a)]
  \item The condition $\PicQ{X}=\NQ{X}$ in the above definition is equivalent to
  $H^1(X,\ox)=0$. Indeed, taking the cohomology of the exponential sequence $0\to \ZZ\to \ox\to \ox^*\to
  0$, one obtains the following exact sequence
   \begin{center}
    \begin{tikzpicture}[>=angle 60,node distance=0.5cm]
     \node (L0) at (0,0) {$0$};
     \node (Pic^0) at (2,0) {$\Pic^0(X)$};
     \node (Pic) at (4.5,0) {$\Pic(X)$};
     \node (NS) at (7,0) {$\NS(X)$};
     \node (R0) at (9,0) {$0$};
     \node (eq) [below of=Pic^0,rotate=90] {$=$};
     \node (subset) [below of=NS,rotate=-90] {$\subseteq$};
     \node [below of=eq] {$H^1(X,\ox)/H^1(X,\ZZ)$};
     \node [below of=subset] {$H^2(X,\ZZ)$};
     \draw [->] (L0) -- (Pic^0);
     \draw [->] (Pic^0) -- (Pic);
     \draw [->] (Pic) -- (NS);
     \draw [->] (NS) -- (R0);
    \end{tikzpicture}
   \end{center}
  where $\NS(X)$ is the group of algebraic equivalence classes of line bundles on $X$. By
  \cite[Remark~1.1.20]{Laz}, a class in $\NS(X)$ is numerically trivial if and only if it
  is torsion, in other words $N^1(X)=\NS(X)_{\text{t.f.}}:=\NS(X)/(\text{torsion})$.
  Thus if we tensor the above sequence with $\QQ$, we get \[
   0\to \Pic^0(X)_\QQ\to \PicQ{X}\to \NQ{X}\to 0, \]
  which shows that $\PicQ{X}= \NQ{X}$ if and only if $H^1(X,\ox)=0$.

  \item Using the same observations in (a), one also sees that suitable conditions
  can imply the even stronger equality $\Pic(X)=N^1(X)$. For example:
   \begin{itemize}
    \item If $\Pic(X)$ is a free abelian group of finite rank (e.g. $X$ is a toric
    variety \cite[\S 3.4, first proposition]{Ful}), then $\Pic(X)=N^1(X)$.
    \item If $X$ is smooth and $H_1(X,\ZZ)=0$ (e.g. Fano variety \cite[Corollary~4.29]{Deb}),
    then $\Pic(X)=N^1(X)$. This is because $H_1(X,\ZZ)=0$ implies $H^1(X,\ox)=0$ by Hodge
    theory, and also implies that $H^2(X,\ZZ)$ is torsion-free by the universal coefficient theorem
    \cite[Corollary~56.4]{Mun}.
   \end{itemize}
 \end{enumerate}
\end{rmk}

\begin{notation} \label{n:V//T}
 For an affine variety $V$ on which an algebraic torus $T$ acts and a character $\chi\colon T\to
 \CC^*$, we will use $V\sslash_{\!\chi} T$ to denote the GIT quotient constructed from
 the $T$-linearized line bundle \[
  \mathcal{O}_V^\chi := \text{the trivial line bundle on $V$, $T$-linearized by $\chi$}.
  \]
 We will also write $V^{\st}_\chi$,
 $V^{\semi}_\chi$, and $V^{\un}_\chi$ to mean the stable, semi-stable, and unstable locus
 of this GIT quotient respectively.
\end{notation}

The next two theorems are among the central results in \cite{HK}. The wording and
notations we use are not exactly the same as the original.

\begin{thm} \label{t:Prop 2.9}
 Let $X$ be a Mori dream space. Let $R=\Cox(X,\LL)$ and let $V$ be the affine variety
 $\Spec R$, with the natural action by the torus $T:=T_{\LL}$ as in
 Definition~\ref{d:Cox}. Let\/ $\chi\in \chi(T)=N^1(X,\LL)$ be a character of\/ $T$
 which corresponds to an ample class in $N^1(X)$. Then $V^{\semi}_\chi$ does not depend on the
 choice of\/ $\chi$, $V\sslash_{\!\chi} T = X$, and the following three properties hold:
 \begin{enumerate}[\upshape (i)]
  \item $V^{\un}_\chi$ has codimension at least $2$ in $V$;
  \item $V^{\semi}_\chi=V^{\st}_\chi$;
  \item Both of the maps \[
    \chi(T)_{\QQ} \rightarrow \Pic^T(V^{\semi}_\chi)_\QQ \leftarrow \Pic(X)_\QQ \]
  are isomorphisms, where the left map sends a character\/ $\nu\in \chi(T)$ to
  $\OO_{V^\semi_\chi}^\nu$, and the right map is the pullback under the quotient
  map $\pi\colon V^{\semi}_\chi \to X$.
 \end{enumerate}
 Moreover, one can choose the basis $\LL$ so that the action of\/ $T$ on $V^\semi_\chi$ is
 free. We will call such basis a preferred basis.
\end{thm}

\begin{proof}
 See the proof of \cite[Proposition~2.9]{HK}.
\end{proof}

\begin{thm} \label{t:Thm 2.3}
 Let\/ $T$ be a torus acting on an affine variety\/ $V$, and let\/ $\chi$ be a character of\/
 $T$. If $X:=V\sslash_{\!\chi} T$ is projective and $\QQ$-factorial, and the conditions
 \textup{(i)--(iii)} of Theorem~\ref{t:Prop 2.9} hold, then $X$ is a Mori dream space.
\end{thm}

\begin{proof}
 See the proof of \cite[Theorem~2.3]{HK}.
\end{proof}

\begin{lemma} \label{l:10}
 Under the same setting as in Theorem~\ref{t:Prop 2.9}, if $\LL$ is a preferred basis,
 then for any line bundle $L$ of the form $L=\LL^\mm$, we have \[
     \pi^* L=\OO_{V^\semi_\chi}^{[L]}  \]
 as $T$-linearized line bundles on $V^\semi_\chi$,
 where $[L]\in N^1(X,\LL)$ is the numerical equivalence class of $L$.
\end{lemma}

\begin{proof}
 Since $T$ acts freely on $V^\semi_\chi$, any $T$-linearized line bundle on
 $V^\semi_\chi$ descends to a (unique) line bundle on $X$; in particular
 $\OO_{V^\semi_\chi}^{[L]}$ descends to a line bundle $M$ on $X$. To identify which
 line bundle $M$ is, we look at the space of $T$-invariant sections of
 $\OO_{V^\semi_\chi}^{[L]}$: on the one hand, via $\pi^*$ we see that it is equal to
 $H^0(X,M)$; on the other hand, we claim that it is also equal to $R_{[L]}$, the space of degree-$[L]$
 homogeneous elements in $R$. Since $R_{[L]}=H^0(X,L)$ by the definition of $R$, we have $M=L$.

 It remains to prove the claim that \[
   \{\text{$T$-invariant sections of $\OO_{V^\semi_\chi}^{[L]}$}\}=R_{[L]}. \]
 Let $\{a_i\in H^0(X,\LL^{\mm_i})\}_{i=1}^\ell$ be a set of regular functions
 on $V$ whose common zero locus is $V^\un_\chi$. A $T$-invariant section $s$ of
 $\OO_{V^\semi_\chi}^{[L]}$ is nothing but a degree-$[L]$ homogeneous regular function on
 $V^\semi_\chi$, so we can represent $s$ as a compatible collection
 $\{b_i/a_i^p\}_{i=1}^\ell$ where $b_i\in H^0(X,L\otimes \LL^{p\mm_i})$, and of course
 the compatibility means \[
     b_i a_j^p = b_j a_i^p, \quad \forall\  i,j\in \{1,\ldots,\ell\}. \]
 As divisors on $X$, we can write $\mathrm{div}(b_i)$ and $\mathrm{div}(a_i^p)$ as
 \begin{align*}
     \mathrm{div}(b_i)&= D_i + B_i, \\
     \mathrm{div}(a_i^p)&= D_i + A_i,
 \end{align*}
 where $D_i$, $B_i$ and $A_i$ are Weil divisors on $X$ such that $B_i$ and $A_i$ have no common
 component. Then the compatibility condition translates to the following equality of
 divisors on $X$: \[
              B_i + A_j = B_j + A_i, \quad \forall\  i,j\in \{1,\ldots,\ell\}. \]
 Since $B_i$ and $A_i$ have no common component, we must have $A_j\geq A_i$, and by symmetry
 $A_i\ge A_j$, thus all the $A_i$'s are the same divisor $A$. But then we must have
 $A=0$, for otherwise if we take a sufficiently large integer $q$ such that $qA$ is
 Cartier and $\OO_X(qA)$ is of the form $\LL^\mm$, then we see that the common zero locus
 of $\{a_i^{pq}\}_{i=1}^\ell$ has codimension one in $V$, contradicting the property~(i) of
 Theorem~\ref{t:Prop 2.9}. Therefore $a_i^p$ divides $b_i$ in $R$ for all $i$, hence the
 section $s$ represented by $\{b_i/a_i^p\}_{i=1}^\ell$ is in $R_{[L]}$.
\end{proof}



\begin{lemma} \label{l:uniqueness}
 Let $X$ be a Mori dream space, and let $\LL$ be a basis of $\PicQ{X}$.
 Suppose that the torus $T:=T_{\LL}$ acts on some normal affine variety $V$, such
 that for some character $\chi\in N^1(X,\LL)$ we have $V\sslash_{\!\chi} T = X$ and the conditions
 \textup{(i)--(iii)} of Theorem~\ref{t:Prop 2.9} hold. Moreover, suppose that
 $\pi^* L=\OO_{V^\semi_\chi}^{[L]}$ for any line bundle $L$ of the form $L=\LL^\mm$,
 where $\pi\colon V^\semi_\chi \to X$ is the quotient map.
 Then $\chi$ corresponds to an ample class
 in $N^1(X)$, and the coordinate ring $R$ of\/ $V$ is equal to $\Cox(X,\LL)$.
\end{lemma}

\begin{proof}
 By \cite[Theorem~8.1]{Dol}, there exists an ample line
 bundle on $X$ whose pullback under the quotient map $\pi\colon V^{\semi}_\chi \to X$ equals some
 tensor power of the $T$-linearized line bundle which was used to construct the GIT
 quotient. It follows from this and the condition (iii) of Theorem~\ref{t:Prop 2.9} that
 $\chi$ corresponds to an ample class.

 Let $L$ be a line bundle on $X$ of the form $\LL^\mm$. Since $\pi^*
 L=\OO_{V^\semi_\chi}^{[L]}$ by assumption, $\pi^*$ induces a natural isomorphism \[
     \{\text{$T$-invariant sections of $\OO_{V^\semi_\chi}^{[L]}$}\}=H^0(X,L). \]
 On the other hand, since $V$ is normal and $V^\un_\chi\subset V$ has codimension $\ge 2$, we have
 \begin{align*}
   \{\text{$T$-invariant sections of $\OO_{V^\semi_\chi}^{[L]}$}\}&=
     \{\text{Degree-$[L]$ homogeneous regular functions on $V^\semi_\chi$}\}  \\
    &=\{\text{Degree-$[L]$ homogeneous regular functions on $V$}\} \\
    &= R_{[L]}.
 \end{align*}
 Hence $R_{[L]}=H^0(X,L)$, and thus  \[
   R=\bigoplus_{L=\LL^\mm} H^0(X,L)= \Cox(X,\LL). \]
\end{proof}

\begin{prop} \label{p:normal}
 With the same setting as in Theorem~\ref{t:Prop 2.9}, if we choose $\LL$ to be a
 preferred basis, then\/ $V=\Spec\Cox(X,\LL)$ is normal.
\end{prop}

\begin{proof}
 Let $\varphi\colon V'\to V$ be the $T$-equivariant normalization of $V$.
 Since $X$ is normal, it follows that $V^{\semi}_\chi$ is also normal (cf. the second paragraph of
 the proof of \cite[Proposition~6.3]{BH}), so $\varphi$ is an isomorphism over
 $V^{\semi}_\chi$.

 We want to show that ${V'}^\un_\chi=\varphi^{-1}(V^\un_\chi)$. The inclusion ${V'}^\un_\chi
 \subset \varphi^{-1}(V^\un_\chi)$ is obvious. For the reverse inclusion, we need to show that if
 $f$ is a homogeneous regular function on $V'$ whose degree is a multiple of $\chi$,
 then $f$ must vanish on $\varphi^{-1}(V^\un_\chi)$. By the definition of normalization,
 $f$ satisfies an integral equation \[
   f^n + a_{n-1} f^{n-1} + \cdots + a_1 f + a_0 = 0, \]
 where the $a_i$'s are homogeneous elements in $R$ whose degrees are multiples of $\chi$.
 So if we plug in a point $p\in \varphi^{-1}(V^\un_\chi)$ into the above equation, we get
 $f(p)^n=0$ since $a_0(p)=\cdots =a_{n-1}(p)=0$, so $f$ vanishes at $p$ as desired.

 Since $\varphi$ is an isomorphism over $V^{\semi}_\chi$ and
 ${V'}^\un_\chi=\varphi^{-1}(V^\un_\chi)$, we thus have ${V'}^\semi_\chi=V^\semi_\chi$,
 so $V'\sslash_{\!\chi} T = X$ and this GIT quotient satisfies the properties \textup{(i)--(iii)}
 of Theorem~\ref{t:Prop 2.9}. Moreover, since $\LL$ is a preferred basis, by
 Lemma~\ref{l:10} $\pi^* L = \OO_{{V'}^\semi_\chi}^{[L]}$ for any line bundle $L$ of the
 form $\LL^\mm$. Hence by Lemma~\ref{l:uniqueness}, the coordinate ring of $V'$ is equal
 to $\Cox(X,\LL)$, namely $V'=V$.
\end{proof}

\section{Proof of Theorem~\ref{t:main}} \label{s:ThmPf}

\begin{proof}[Proof of Theorem~\ref{t:main}]
 We need to verify that $Y\subset X$ satisfies the three requirements in
 Definition~\ref{d:1}. First, by the Lefschetz hyperplane theorem
 \cite[Example~3.1.24 and Example~3.1.25]{Laz},
 the restriction map determines canonical isomorphisms $N^1(X)=N^1(Y)$ and
 $\Pic(X)=\Pic(Y)$. In particular, the second requirement in Definition~\ref{d:1} is
 satisfied, and $\PicQ{Y}=\NQ{Y}$ since $\PicQ{X}=\NQ{X}$.

 We pick the basis $\LL$ to contain the line bundle $\ox(Y)$. Let $R=\Cox(X,\LL)$, and let
 $s\in R$ be the unique (up to constant multiples) section of
 $\ox(Y)$ whose zero locus is $Y$. Since $Y$ is irreducible, the ideal $sR\subset R$ is a prime
 ideal, so it defines an irreducible subvariety $W\subset V$. By Theorem~\ref{t:Prop
 2.9}, $V^{\un}_\chi$ does not depend on the choice of ample $\chi\in N^1(X,\LL)$, and since
 $Y$ is ample, we have $V^{\un}_\chi\subset W$. Therefore \[
  W^{\un}_\chi = V^{\un}_\chi \cap W = V^{\un}_\chi \]
 has codimension $\geq 2$ in $W$ due to the condition~\eqref{e:*}. Also by
 Theorem~\ref{t:Prop 2.9},
 $V^{\semi}_\chi=V^{\st}_\chi$ and $V\sslash_{\!\chi} T = V^{\st}_\chi/T = X$ is a good
 geometric quotient, so it follows that $W^{\semi}_\chi=W^{\st}_\chi$ and
 $W\sslash_{\!\chi} T = W^{\st}_\chi/T = Y$ is a good geometric quotient. Thus we see
 that the GIT quotient $W\sslash_{\!\chi} T = Y$ satisfies the conditions (i) and (ii) of
 Theorem~\ref{t:Prop 2.9}, and we claim that it satisfies the condition (iii) as well:
 the map $\Pic^T(W^{\semi}_\chi)_\QQ \leftarrow \Pic(Y)_\QQ$ is an isomorphism by
 Kempf's descent lemma \cite[Th\'{e}or\`{e}me~2.3]{DN}, and thus the map
 $\chi(T)_{\QQ} \rightarrow \Pic^T(W^{\semi}_\chi)_\QQ$ is also an isomorphism since
 $\chi(T)_{\QQ}=\PicQ{X}=\PicQ{Y}$. Therefore $Y$
 is a Mori dream space by Theorem~\ref{t:Thm 2.3}, which verifies the first requirement in
 Definition~\ref{d:1}.

 To verify that $Y$ respects the chamber structure, we will use the fact that the Mori
 chambers coincide with the GIT chambers \cite[Theorem~2.3]{HK}. Suppose $\chi_1$,
 $\chi_2\in N^1(X)$ are in the interior of the same Mori chamber of $X$. Then
 $V^{\un}_{\chi_1}=V^{\un}_{\chi_2}$, so \[
  W^{\un}_{\chi_1} = V^{\un}_{\chi_1} \cap W = V^{\un}_{\chi_2} \cap W =
  W^{\un}_{\chi_2}. \]
 Hence $\chi_1$ and $\chi_2$ are also in the same Mori chamber of $Y$.

 To show that $\Nef(X)=\Nef(Y)$, suppose on the contrary that $\Nef(Y)\supsetneqq\Nef(X)$.
 Then there exists $\nu\in N^1(X)$ which is ample on $Y$ and lies in the interior of some
 Mori chamber of
 $X$ not equal to $\Nef(X)$. Since $\nu$ and $\chi$ are both ample on $Y$, $W^\un_\nu = W^\un_\chi$.
 Recalling that $W^{\un}_\chi = V^{\un}_\chi$, we thus have \[
   V^{\un}_\nu \cap W = W^{\un}_\nu = W^{\un}_\chi = V^{\un}_\chi, \]
 so $V^{\un}_\nu \supset V^{\un}_\chi$. Since $V^{\semi}_\chi \sslash T = X$ is a good
 geometric quotient, if $V^{\un}_\nu \supsetneqq V^{\un}_\chi$ then $V^{\semi}_\nu \sslash
 T$ would be an open subset of $X$ and hence not projective, a contradiction. So $V^{\un}_\nu =
 V^{\un}_\chi$, and hence $V^{\semi}_\nu \sslash T = V^{\semi}_\chi \sslash T = X$, which
 means $\nu$ and $\chi$ are both in $\Nef(X)$.
\end{proof}

\section{Proof of Proposition~\ref{p:4} and \ref{p:5}} \label{s:PropPf}

\begin{proof}[Proof of Proposition~\ref{p:4}]
 (a) If the Picard number of $X$ is one, then we can pick a very ample line bundle $L$ on
 $X$ as a basis of $\PicQ{X}$ so that $R:=\Cox(X,L)$ is precisely the
 homogeneous coordinate ring of $X$ embedded into some projective space by the linear series $|L|$.
 Hence $V:=\Spec R$ is the corresponding affine cone over $X$, $V^\un_{[L]}$ is the origin, and
 $T=\CC^*$. The codimension of $V^\un_{[L]}$ in $V$ is thus equal to $\dim V=\dim X+1\geq 3$.

 (b) By Remark~\ref{r:Pic=N}~(a) we have $H^1(X_1,\OO_{X_1})=0$, which implies
 $\Pic(X_1\times X_2)=\Pic(X_1)\times \Pic(X_2)$ \cite[Chapter~III
 Exercise~12.6]{Har}, and hence also $N^1(X_1\times X_2)=N^1(X_1)\times N^1(X_2)$.
 Let $\LL_i$ be a basis of $\PicQ{X_i}$, $R_i=\Cox(X_i,\LL_i)$, $V_i=\Spec R_i$,
 $T_i=\Hom (N^1(X_i,\LL_i),\CC^*)$, and
 $p_i\colon X_1\times X_2 \to X_i$ be the projection map for $i=1,2$. Let $\LL$ be the
 basis of $\PicQ{X_1\times X_2}$ which consists of $p_1^*\LL_1$ and $p_2^* \LL_2$. Then
 by the K\"{u}nneth formula we have \[
        R:=\Cox(X_1\times X_2,\LL)=R_1\otimes R_2. \]
 Thus $V:=\Spec R = V_1\times V_2$ and $T:= \Hom (N^1(X_1\times X_2,\LL),\CC^*)=T_1\times
 T_2$. Furthermore, since $\chi= p_1^*\,\chi_1 + p_2^*\,\chi_2$ is ample if and only if
 $\chi_i\in N^1(X_i)$ are both ample for $i=1,2$, we thus have \[
   V^\un_\chi = p_1^{-1} {V_1^\un}_{\!\!\!\!\!\chi_1} \,\cup\, p_2^{-1} {V_2^\un}_{\!\!\!\!\!\chi_2}, \]
 hence the result follows.

 (c) We will in fact show that if the line bundles $L_i$'s are all of the form
 $\LL^{\mm_i}$ where $\LL$ is a \emph{preferred basis}, then
 $\widetilde{X}:=\PP(\bigoplus_{i=1}^k L_i)$ is a Mori dream space which has small
 unstable locus.
 Let $R=\Cox(X,\LL)$, $V=\Spec R$, $T=\Hom (N^1(X,\LL),\CC^*)$, and let $\chi\in N^1(X,\LL)$ be a
 character of $T$ which corresponds to a sufficiently ample class, to the extent that \[
   \chi + \sum_{i=1}^k d_i [L_i]  \]
 is ample for all $d_i\ge 0$ and $\sum_{i=1}^k d_i=1$. Let $\widetilde{\LL}$ be the basis
 of $\PicQ{\widetilde{X}}$ consisting of $p^* \LL$ and $\OO_{\widetilde{X}}(1)$, where
 $p:\widetilde{X}\to X$ is the projection map. Then \[
  \widetilde{T}:= \Hom (N^1(\widetilde{X},\widetilde{\LL}),\CC^*)= T\times T_1, \]
 where $T_1=\CC^*$ is the one-dimensional torus for which $[\OO_{\widetilde{X}}(1)]$ generates
 the group of characters.

 By Lemma~\ref{l:10}, we have $\pi^* L_i = \OO_{V^\semi_\chi}^{[L_i]}$, $i=1,\ldots,k$.
 Thus we consider the normal affine variety  \[
  \widetilde{V}:= V \times \overbrace{\CC \times \CC \times \cdots \times \CC}^{\text{$k$
  copies}} \]
 (the normality of $V$ follows from Proposition~\ref{p:normal}) with the following
 $\widetilde{T}$-action: given any
 $(v,a_1,\ldots,a_k)\in \widetilde{V}$ and $(t,c)\in T\times \CC^*=\widetilde{T}$, define  \[
  (t,c)\cdot (v,a_1,\ldots,a_k) = (t\cdot v,\<[L_1],t\>^{-1}a_1c,\ldots,\<[L_k],t\>^{-1}a_kc), \]
 where $\<[L_i],t\>$ denotes the natural pairing. Let $\widetilde{\chi}$ be the character
 $(\chi,[\OO_{\widetilde{X}}(1)])$ of\/ $\widetilde{T}$. We claim that  \[
  \widetilde{V}^\un_{\widetilde{\chi}}= p_1^{-1} V^\un_\chi\,\cup\,p_2^{-1}
  \{(0,\ldots,0)\}, \]
 where $p_1\colon \widetilde{V}\to V$ and $p_2\colon \widetilde{V}\to \overbrace{\CC \times \cdots \times \CC}^{\text{$k$
 copies}}$ are the projection maps. To see this, let $x_1,\ldots,x_k$ be the coordinate
 functions of $\overbrace{\CC \times \cdots \times \CC}^{\text{$k$ copies}}$, and let
 $f\in R_\nu$ be a homogeneous degree-$\nu$ regular function on $V$ for some $\nu\in
 N^1(X,\LL)$. If $d_1,\ldots,d_k$ are nonnegative integers which sum up to $d$, then
 $f x_1^{d_1}\cdots x_k^{d_k}$ is a regular function on $\widetilde{V}$ which is
 homogeneous of degree \[
            (\nu-\sum_{i=1}^k d_i[L_i],d[\OO_{\widetilde{X}}(1)]),   \]
 hence $\widetilde{V}^\un_{\widetilde{\chi}}$ is the common zero locus of all such
 functions $f x_1^{d_1}\cdots x_k^{d_k}$ for which  \[
   \nu = d\chi+\sum_{i=1}^k d_i[L_i]. \]
 Note that such $\nu$ corresponds to an ample class thanks to our choice of $\chi$ in the
 beginning, and from this it follows easily that \[
  \widetilde{V}^\un_{\widetilde{\chi}}= p_1^{-1} V^\un_\chi\,\cup\,p_2^{-1}
  \{(0,\ldots,0)\}. \]
 Hence $\widetilde{V}^\un_{\widetilde{\chi}}$ has codimension $\ge 3$ in $\widetilde{V}$,
 and $\widetilde{V}\sslash_{\!\widetilde{\chi}} \widetilde{T}=\PP(\bigoplus_{i=1}^k
 L_i)$, so in particular this GIT quotient satisfies the property~(i) in Theorem~\ref{t:Prop 2.9}.
 The action of $T$ on $V^\semi_\chi$ is free since $\LL$ is a preferred basis, hence the
 action of $\widetilde{T}$ on $\widetilde{V}^\semi_{\widetilde{\chi}}$ is also free, so
 in particular the GIT quotient $\widetilde{V}\sslash_{\!\widetilde{\chi}} \widetilde{T}$
 satisfies the property~(ii) in Theorem~\ref{t:Prop 2.9}, and the right map in the
 property~(iii) is an isomorphism even before tensoring with $\QQ$; to show the left map is also
 an isomorphism, we will show that $ \OO_{\widetilde{V}^\semi_{\widetilde{\chi}}}^{[\widetilde{L}]}
 =\widetilde{\pi}^* \widetilde{L}$, where
 $\widetilde{\pi}\colon \widetilde{V}^\semi_{\widetilde{\chi}}\to \widetilde{X}$ is the
 quotient map, and $\widetilde{L}$ is any line bundle on $\widetilde{X}$ of the form
 $\widetilde{\LL}^{\widetilde{\mm}}$. Since we already know that the pullback map
 $\Pic^{\widetilde{T}}(\widetilde{V}^{\semi}_{\widetilde{\chi}}) \leftarrow
 \Pic(\widetilde{X})$ is an isomorphism, the line bundle
 $\OO_{\widetilde{V}^\semi_{\widetilde{\chi}}}^{[\widetilde{L}]}$ descends to some line
 bundle $\widetilde{M}$ on $\widetilde{X}$ in any case, so we just need to identify which
 line bundle $\widetilde{M}$ is. To do this we look at the space of $\widetilde{T}$-invariant
 sections of $\OO_{\widetilde{V}^\semi_{\widetilde{\chi}}}^{[\widetilde{L}]}$: on the one
 hand it is equal to $H^0(\widetilde{X},\widetilde{M})$ via $\widetilde{\pi}^*$; on the
 other hand, it is equal to the space of homogeneous degree-$[\widetilde{L}]$ regular
 functions on $\widetilde{V}$ (since $\widetilde{V}$ is normal), and if
 $\widetilde{L}=p^* L \otimes \OO_{\widetilde{X}}(d)$, then such regular functions are
 precisely
 linear combinations of functions of the form $f x_1^{d_1}\cdots x_k^{d_k}$ where the
 $d_i$'s are nonnegative integers summing up to $d$ and $f$ is a homogeneous regular function
 on $V$ of degree $[L]+\sum_{i=1}^k d_i[L_i]$. In other words, \[
  \{\text{$\widetilde{T}$-invariant sections of
  $\OO_{\widetilde{V}^\semi_{\widetilde{\chi}}}^{[\widetilde{L}]}$}\}
      =\bigoplus_{\substack{d_1,\ldots,d_k \ge 0\\ d_1+\cdots+d_k=d}}
                   H^0(X,L\otimes L_1^{\otimes d_1}\otimes\cdots\otimes L_k^{\otimes d_k}). \]
 But the right-hand side is exactly $H^0(\widetilde{X},\widetilde{L})$: indeed since
 $\widetilde{L}=p^* L \otimes \OO_{\widetilde{X}}(d)$, by the projection formula
 \begin{align*}
  H^0(\widetilde{X},\widetilde{L})&=H^0(X,p_*\widetilde{L})
   =H^0(X,L\otimes p_*\OO_{\widetilde{X}}(d))=H^0(X,L\otimes S^d(\bigoplus_{i=1}^k L_i)) \\
   &=\bigoplus_{\substack{d_1,\ldots,d_k \ge 0\\ d_1+\cdots+d_k=d}}
                   H^0(X,L\otimes L_1^{\otimes d_1}\otimes\cdots\otimes L_k^{\otimes d_k}).
 \end{align*}
 Hence the line bundle on $\widetilde{X}$ which
 $\OO_{\widetilde{V}^\semi_{\widetilde{\chi}}}^{[\widetilde{L}]}$ descends to must be
 $\widetilde{L}$.

 Now we can use Theorem~\ref{t:Thm 2.3} to obtain that $\widetilde{X}$ is a Mori dream space,
 and then use Lemma~\ref{l:uniqueness} to conclude that the affine variety
 $\widetilde{V}$ is indeed $\Spec \Cox(\widetilde{X},\widetilde{\LL})$, thus completing the
 proof.
\end{proof}

\begin{proof}[Proof of Proposition~\ref{p:5}]
 It was shown in \cite{Cox} that in the toric case, the Cox ring $R$ and the unstable
 locus $V^\un_\chi$ have the following explicit description. Let $\Delta(1)$ be the
 set of all one-dimensional cones of $\Delta$. For each $\rho\in\Delta(1)$,
 introduce a variable $x_\rho$. Then $R$ is the polynomial ring \[
      R=\CC[x_\rho\colon \rho\in\Delta(1)], \]
 and the unstable locus is the zero locus of an ideal $I\subset R$ generated by squarefree monomials.
 To give these generators of $I$, let us introduce some notations. For a subset $\AA
 \subset \Delta(1)$, we denote the monomial $\prod_{\rho\in\AA}x_\rho$ as $\xx^\AA$, and
 we write $\widehat{\AA}$ for the complement of $\AA$ in $\Delta(1)$.
 For a cone $\sigma\in \Delta$, we let
 $\sigma(1)=\{\rho\in\Delta(1)\colon \rho \text{ is a face of }\sigma\}$. Then \[
  I=\< \xx^{\widehat{\sigma(1)}}\colon \sigma\in\Delta \> \subset R. \]

 A squarefree monomial ideal such as $I$ can be very well understood by associating to it a
 simplicial complex, sometimes called the Stanley-Reisner complex. The Stanley-Reisner
 complex of $I$ is by definition an abstract simplicial complex $\Sigma$ on the vertex
 set $\Delta(1)$, whose faces are those subsets $\AA \subset \Delta(1)$ such that $\xx^\AA \notin
 I$. By \cite[Theorem~1.7]{MS}, the prime decomposition of $I$ is given by \[
  I= \bigcap_{\AA\in\Sigma} \< x_\rho \colon \rho\in \widehat{\AA} \>. \]
 Hence
 \begin{align*}
  \text{The zero locus of } &I \text{ has codimension $\geq 3$} \iff \text{$|\widehat{\AA}|\geq
   3$ for all $\AA\in\Sigma$} \\
   &\iff \text{If $\AA\subset\Delta(1)$ such that $|\widehat{\AA}|\leq 2$, then $\xx^\AA\in
   I$} \\
   &\iff \text{If $\AA\subset\Delta(1)$ such that $|\AA|\leq 2$, then $\xx^{\widehat{\AA}}\in
   I$} \\
   &\iff \text{If $\AA\subset\Delta(1)$ and $|\AA|\leq 2$, then $\AA=\sigma(1)$ for some $\sigma\in\Delta$.}
 \end{align*}
\end{proof}


\end{document}